\documentclass[12pt,reqno,english,a4]{amsart}

\usepackage{amsmath,amsthm,amsfonts,amssymb,amscd}
\usepackage[latin1]{inputenc}
\usepackage{psfrag}
\usepackage{epsfig}
\usepackage{a4wide}

\usepackage{hyperref} 
\theoremstyle{plain}
\newtheorem{theorem}{Theorem}[section]

\newtheorem{lemma}[theorem]{Lemma}
\newtheorem{claim}[theorem]{Claim}
\newtheorem{example}[theorem]{Example}

\newtheorem{definition}[theorem]{Definition}

\numberwithin{equation}{section}

\newcommand{\ep} {\epsilon}

\renewcommand{\th} {\theta}

\newcommand{\ra}{\rightarrow}

\newcommand{\om} {\omega}       

\def \RR {{\mathbb R}}
\def \ZZ {{\mathbb Z}}
\def \NN {{\mathbb N}}

\def \TT {{\mathbb T}}
\def \SS {{\mathbb S}}

\newcommand{\cO}{{\mathcal O}}

\begin{document}

\title[Robustly transitive actions]
      {Robustly transitive actions of $\RR^2$ on  compact three manifolds}

\author{Ali Tahzibi}\thanks{The authors would like to thank Fapesp
(Fundação de Amparo a Pesquisa de Estado de São Paulo) for
financial support (projeto temático 05/03107-9.)}

\author{Carlos Maquera}\thanks{}



\date{\today}

\address{Ali Tahzibi and Carlos Maquera\\ Universidade de
S{\~a}o Paulo - S{\~a}o Carlos \\Instituto de ci{\^e}ncias
matem{\'a}ticas e de Computa\c{c}{\~a}o\\
Av. do Trabalhador S{\~a}o-Carlense 400 \\
13560-970 S{\~a}o Carlos, SP\\
Brazil}

 \email{tahzibi@icmc.usp.br}
 \email{cmaquera@icmc.usp.br}

 \maketitle

 \begin{abstract}
 In this paper, we define $C^1$-robust transitivity for  actions of
 $\RR^2$ on closed connected orientable  manifolds.  We prove that if
 the ambient manifold is three dimensional and the dense orbit of a robustly transitive action is
 not  planar, then it is
  ``degenerate" and  the action is defined by an Anosov flow.
 \end{abstract}

 \medskip
 \medskip

 \thispagestyle{empty}


 \section{Introduction}
 In some recent works in the theory of dynamical systems  robust transitivity of diffeomorphisms and
 flows has been investigated. Weak forms of hyperbolicity has been shown to be necessary conditions
 for robust transitivity of flows and diffeomorphisms of compact manifolds. Bonatti-Díaz-Pujals \cite{BDP} proved that
 $C^1$-robustly transitive diffeomorphisms admits dominated splittings. Previous to this work,
 Díaz-Pujals-Ures \cite{DPU99} had proved
 that robustly transitive diffeomorphism on three dimensional manifolds are partiall hyperbolic.
 For $C^1$-flows, there are also results which imply that some weak form of
 hyperbolicity is necessary to obtain robust transitivity.
  See for example a result of Vivier \cite{vivier} about robustly transitive flows
 on any dimension and a result of Doering \cite{Do87} in three dimensional case.

    By using Kupka-Smale theorem in
 one dimensional case, one deduce that there does not exist any
 robustly transitive  diffeomorphism on one dimensional manifolds. Also, by
 a result of Peixoto \cite{peixoto}, we know that the Morse-Smale flows form a
 dense subset of the set of $C^1$-flows on any surface. This result
 implies that robustly transitive flows may exist only on manifolds with
 dimension higher than two.

 If we consider the diffeomorphisms or flows defined on a manifold
 as the action of $\ZZ, \, \RR$ on it, a naturally question arises:
``what about robustly transitive actions of higher dimensional
groups?"

In this paper, we begin the study of robustly transitive actions
of $\RR^2$ by giving some examples of these actions and proving
that in three dimensional manifolds the only robustly transitive
actions of $\RR^2$ (We do not consider the case when all orbits
are planar) are defined by robustly transitive flows (see theorem
\ref{teo:main}). On the other hand we know that by a result of
Doering, robustly transitive flows are Anosov flows.

Let $N$ denote a closed connected orientable three manifold and
$\varphi:\mathbb{R}^{2}\times N\to N$\, be a \,$C^{r}$-action.
 For each \,$w \in \mathbb{R}^2 \setminus \{0 \}, \varphi $
 induces a  \,$C^{r}$-flow \,$ ( \varphi_{w}^{t})_{t\in\mathbb{R}} $\, given by
  \,$\varphi_{w}^{t}(p)=\varphi (tw,p)$\,  and its corresponding \,$C^{r-1 }$-vector field
  \,$X_{w}$\, is defined by \,$X_{w}(p)= D_{1}\varphi(0,p)\cdot w.$\,
 If \,$\{w_{1},w_{2}\}$\, is a base of \,$\mathbb{R}^{2},$\, the associated vector fields
 \,$X_{w_{1}},X_{w_{2}}$\, satisfy the commutativity condition \,$[X_{w_{1}},X_{w_{2}}]=0$\, and determine completely  the action
 \,$\varphi.$\, They are called \emph{infinitesimal generators} of \,$\varphi.$ This condition of commutativity
 between two vector fields is a necessary and sufficient condition for them to be generators of an action.
 \,$X_{(1,0)}$\, and \,$X_{(0,1)}$\, are called the \emph{canonical infinitesimal generators}.

 Denote by $A^{r}(\mathbb{R}^{2},N)\,\ 1\leq r\leq \omega \,$
 the set of
 actions of \,$\mathbb{R}^{2}$\, on \,$N$\, whose infinitesimal generators are of
 class \,$C^{r}.$\,
 Given two actions \,$\{ \varphi; X_{(1,0)},X_{(0,1)}\}$\, and \,$\{ \psi; Y_{(1,0)},Y_{(0,1)}\}$\,
 define \,$$d_{(1,1)}(\varphi,\psi)= \max\{\|X_{(1,0)} - Y_{(1,0)}\|_{1},\|X_{(0,1)} - Y_{(0,1)}\|_{1}\}.$$\
 With this distance \,$A^{r}(\mathbb{R}^{2},N)$\, is a metric space and
 its corresponding topology is called the \emph{\,$C^{(1,1)}$-topology.}
 Note that this topology is finer than the \,$C^{2}$-topology  and
 coarser than the \,$C^{1}$-topology.
For any action $\phi \in A^r(\RR^2, M),$ $\mathcal{O}_p :=  \{
\phi(\om, p), \om \in \RR^2 \}$ is called the orbit of $p \in M.$
 The orbit is called
\textit{singular} if its dimension is less than two.
\begin{definition}
 { \rm An action $\varphi \in A^1(\mathbb{R}^{2},N)$ is called
 \textit{transitive} if it admits a dense orbit in $N.$
 $\varphi$ is \textit{robustly transitive} if it admits a $C^{(1,1)}$ neighborhood
 $\mathcal{U}$ of transitive actions.}
 \end{definition}

Our main result is the following theorem.
 \begin{theorem}[Main Theorem]
 \label{teo:main}
 Let $N$ be a closed orientable $3$-manifold. Assume that
 $\varphi \in A^{2}(\mathbb{R}^{2},N)$ is robustly transitive with a
 dense orbit which is not homeomorphic to $\RR^2$. Then,
 $\varphi$ is defined by an Anosov flow.
 \end{theorem}

We mention that the hypotheses about the topological type of the
dense orbit is important to our result. By this hypotheses, the
dense orbit is cylindrical or homeomorphic to $\RR.$ However, we
conjecture that the same result is true without this hypotheses.\\
We would like to thank C. Bonatti for mentioning  us that,
Rosenberg had left the stability problem of the action with all
leaves planar (which is the case we are avoiding here) as an open
problem.

To prove the above theorem, we study the topological type of the
obits of a robustly transitive action. In general one can have three
different  topological type for non singular (two dimensional)
orbits of an action of $\RR^2.$ But in the context of transitive
action, we show that the topological types are more restricted. We
show that for the purpose of proving our theorem, we can suppose
that the action has a dense cylindrical orbit. Then, we apply a
closing lemma for actions (proves by Roussarie for locally free
actions) and prove that $C^{(1, 1)}$-close to a robustly transitive
action, one can find an action with many compact leaves. The
existence of sufficiently large amount of these tori will contradict
the existence of a dense leaf and consequently robust transitivity
of the initial action.

\section{Examples and basic Results}

Let us give some examples of robustly transitive actions of
$\RR^2$.
\begin{example}
{\rm Firstly we construct a singular (defined by flow) example of
a robustly transitive
 action. Consider a robustly transitive Anosov vector field $X$ defined
 in $\TT^3$ and let $\phi \in A(\RR^2, \TT^3)$ be the action defined
 by $X_1 := X$ and $X_2 := c X (c \in \RR).$

 It is obvious that $[X_1,X_2] = 0$ and so they define an action of $\RR^2$ in
 $\TT^3.$ Clearly, all orbits of  this action are singular.
  We claim that $\phi$ is robustly
 transitive. Indeed, suppose $\psi \in A(\RR^2, \TT^3)$ any $C^{(1,1)}$ perturbation of
 $\phi$. By the definition of $C^{(1,1)}-$topology in $A(\RR^2, \TT^3)$ we
 conclude that $\tilde{\phi}$ is defined by $\tilde{X_1},
 \tilde{X_2}$ such that $[\tilde{X_1}, \tilde{X_2}] = 0$ and $\tilde{X_1}$ is $C^1$- close to
 $X_1.$ So, $\tilde{X}_2$ is also an Anosov flow.
  By a result of Kato-Morimoto \cite{kamo73} the
 centralizer of an Anosov vector field is trivial. This means that $\tilde{X_2} =
 f
 \tilde{X_1}$ ($f$ is a first integral.) and consequently $\psi$ is also defined by a
 transitive flow. On the other hand, let $Y$ be any flow $C^1$-close to $X$, then $Y$ and
 $cY$ defines an action which is $C^{(1,1)}-$close to $\phi$ and consequently transitive.
 This implies that $Y$ is a transitive flow and so, $X$ is $C^1$-robustly
 transitive. So, we showed that $\phi$ is defined by a robustly
 transitive flow.}
\end{example}

Let denote by $X^t$ the flow of a vector field $X.$

\begin{example}
 {\rm Let $N$ be a three manifold supporting a robustly transitive Anosov flow.
 In the second example, we construct a robustly transitive action in
 $M^4 = N \times \SS^1$ which is not defined by a flow. Consider the coordinate system
 $(x, \th)$ in $M^4, x \in N, \th \in \SS^1.$ In what follows, for a real function
 $a(x, \th)$,  by  $ a(x, \th)\frac{\partial}{\partial x} $ we mean
 $ a_1 \frac{\partial}{\partial x_1} + a_2 \frac{\partial}{\partial x_2}
 + a_3 \frac{\partial}{\partial x_3}$ where $x_1,x_2, x_3$ are coordinates in $N.$

 Let $\phi \in A(\RR^2, M^4)$ be defined by $X_1$ and $X_2$  such that
 $X_1 = a(x, \th) \frac{\partial}{\partial x}$ is a robustly transitive Anosov flow
 in $N$ and $X_2:= \frac{\partial}{\partial \th}.$
 We claim that $\phi$ is robustly transitive.

 Consider a $C^{(1,1)}$-perturbation $\psi$ of the initial action
 $\phi.$ It is generated by two vector fields $Y_1$ and
 $Y_2$ which are respectively  $C^1$ close to $X_1$ and
 $X_2.$

 Let $N_0 := \{ (x, 0) : x \in N\}.$ By transversality of $X_2$ to $N_0$ and
 closedness of $X_2$ and $Y_2$ we conclude that $Y_2$ is also transverse to $N_0.$

 In our coordinate systems
 $$
 Y_1 = \tilde{a}(x, \th) \frac{\partial}{\partial x} + b(x, \th) \frac{\partial}{\partial
 \th}.
 $$
 $$
 \tilde{X_2} = c(x, \th) \frac{\partial}{\partial x} + d(x, \th) \frac{\partial}{\partial
 \th}.
 $$
where $b$ and $c$ are close to zero in $C^1$-topology, $a$ and
$\tilde{a}$ are close in each coordinates and $d$ is close to
constant $1$. We define
$$ \Pi (Y_1) = d Y_1 - b Y_2 = (\tilde{a} d - b c) \frac{\partial}{\partial x}.$$

Observe that $\Pi(Y_1)$ is a $C^1$-vector field close to $X_1$ and
consequently it is transitive. The intersection of the orbits of
$\psi$ with $N_0$ coincide with the orbits of $ \Pi (Y_1).$
 Let $x_0 \in N_0$ with a  $\Pi
(Y_1)$ dense orbit. We claim that the orbit of $\psi$ passing
through $(x_0, 0)$ is dense in $M^4.$ To see this, just observe that
$N_0$ is a global transverse manifold for $Y_2.$ Let $U \in M^4$ an
open set and $V = \bigcup_{t \in \RR} Y_2^t(U) \cap N_0.$ Then, $V$
is an open subset of $N_0$ and by transitivity of $\Pi
(\tilde{X_1})$ there exists $x \in N_0$ such that  for some $t \in
\RR$, $Y_1^t (x) \in V$ and consequently for some $s \in \RR, Y_2^s
(Y_1^t (x)) \in U $ and this finishes the proof. }
\end{example}
Here we outline some basic results about the action of $\RR^2$
which will be used in our proof of the main result. Recall that,
for any action $\phi \in A^1(\RR^2, M),$ $\mathcal{O}_p :=  \{
\phi(\om, p), \om \in \RR^2 \}$ is the orbit of $p \in M$ and $G_p
:= \{ \om \in \RR^2 : \phi(\om, p) = p  \}$ is called the isotropy
group of $p.$ One of the simple results is the following.

\begin{lemma}  Suppose that $\cO (q) $ is accumulated by $\cO (p)$
then $G_p \subseteq G_q.$
 \end{lemma}
 \begin{proof}
 To prove, just observe that for $\om \in G_p$, by definition of action and isotropy group
  we have $\phi(\om, \phi(\eta, p) ) = \phi(\eta, p).$ So, by
  continuity of $\phi$ we conculde that, if $z$ is an accumulation
  point of $\cO (p)$ then $\phi (\om , z) = z$ and consequently  we have $\om \in G_q.$
 \end{proof}

  Using the above lemma we can show that all the dense orbits of $\varphi$ have the same topological type.
  In the setting of Theorem \ref{teo:main} all the dense orbits are
  either
  line or cylinder. In fact we point out that the existence of a
  dense line prohibits the existence of any (not necessarily dense)
  cylinder. This is just by the continuity of action. Consequently,
  if the dense orbit for a transitive action is homeomorphic to $\RR$
  then it is given by a transitive flow.

 \begin{claim} \label{seila} Any two dense orbits are homeomorphic.
 \end{claim}
 \begin{proof} To prove the
claim let $\mathcal{O}_p, \mathcal{O}_q$ two dense orbits with
isotropy groups $G_p, G_q.$ The density of $\mathcal{O}_p$ implies
that $G_p \subset G_q$. Similarly the density of $\mathcal{O}_q$
implies that $G_q \subset G_p$ and consequently $G_p = G_q$ and the
orbits are homeomorphic.
\end{proof}

 \begin{lemma} \label{seila2}
  If $\phi \in A^1(\RR^2, M)$ has a dense cylinder orbit, then any
  two dimensional orbit is either torus or cylinder.
 \end{lemma}

\begin{proof}
 Suppose that $\mathcal{O}_p$ is a dense cylinder orbit. The isotropy
 group of $p$ is $ \ZZ u $ for some $ 0 \neq u  \in \RR^2.$ Let
 $Y$ be the vector field whose flow corresponding is $Y^t = \phi(tu,\cdot)$. It is clear
 that $Y^1(p) = p.$ Let $X$ be the vector field whose flow corresponding is
 $X^t = \phi(tv,\cdot)$, where $v$ is any linearly independent to $u$.
 Then $X$ and $Y$ commute and consequently every point on
 $\mathcal{O}_p$ is periodic with period one for $Y.$ Now,
 using denseness of $\mathcal{O}_p$ we conclude that any point of
 the manifold is a periodic point for $Y$ which finishes the proof
 of the lemma.
\end{proof}

\section{Closing lemma and proof of the main result}
  First of all let us recall the closing lemma of Pugh (\cite[Theorem 6.1]{pu67a})  for the flows in a
  two
  dimensional manifold.
  \begin{theorem}  \label{t.closing}
  Let $X \in \mathcal{X}(M^{2})$ have a nontrivial recurrent
  tranjectory through $p^{*} \in M,$ let $U$ be a neighbourhood of
  $p^{*}$ and $\ep > 0$ be given. Then, there exists $Z \in \mathcal{X}
  (M)$ such that:
  \begin{enumerate}
   \item $X - Z$ vanishes om $ M \setminus U,$
   \item the $C^1$-size of $X- Z$ is less than $\ep$ respecting
   the $U$-coordinates,
   \item $Z$ has a closed orbit through $p^*.$
  \end{enumerate}

  \end{theorem}
In \cite{Ro1}, Roussarie and Weil proved a closing lemma for the action of $\mathbb{R}^2$ on three manifolds.
More precisely one of their results is the following:
\begin{theorem}
\label{roussarie} Let $N$ be an orientable compact closed $C^r$ $(r
> 2)$, $3$-manifold and $\varphi$ a
locally free $C^r$-action. If all orbits of $\varphi$ are not
planar, then there is a locally free action $\varphi_1 \in A^{r}
(\mathbb{R}^{2} ,N)$ with a compact orbit and $C^1-$close to
$\varphi.$

\end{theorem}

To prove the above theorem, the authors firstly observe that either
$\varphi$ has a compact orbit or all the orbits are dense. In the
latter case just take $\varphi = \varphi_1$. In the former case, the
denseness of all orbits is a corollary of a result of Sacksteder
\cite{sack}
 about the minimal sets of $\mathbb{R}^{n-1}$ actions on
$n-$manifolds. The result of Saksteder states that there is no
exceptional minimal set for locally free actions. Using the
denseness of a cylinder, one can show that all other orbits are
cylindrical. In this setting (all the orbits are cylindrical) the
proof of Pugh closing lemma for flows on surfaces can be carried
on to prove that $\varphi$ can be perturbed to give a compact
orbit.

 Let us mention that the above theorem is not the main result of Roussarie and Weil´s paper.
 In fact, their paper is mainly dedicated to the proof of the following theorem \cite[Theorem 2
 (1)]{Ro1}.

\begin{theorem} Let $\phi$ be a $C^r$-action. For all non-planar and
recurrent orbit $\Lambda$ e for all $\ep
> 0$ there exists a submanifold diffeomorphic to $\TT^2$,
$\ep$-close to $\Lambda$ such that the plane field tangent to this
submanifold can be extended to a plane field $C^1$ near to plane
field correponding to $\phi.$ \end{theorem} The main issue in this
result is to find a nearby torus to the recurrent leaf.
 Here we
have a general action which can have singularities. However, we
suppose that there exists a dense cylinder. The idea is using
again the closing lemma of Pugh.

\begin{theorem}
\label{teo:closingLem } Let $N$ be an orientable compact closed,
$3$-manifold and $\varphi \in A^{1} ( \mathbb{R}^{2} , N)$. If
there exists a dense orbit $\mathcal{L}$ of $\varphi$ homeomorphic
to $\SS^1 \times \mathbb{R}$, then there is an action $\varphi_1
\in A^{1} (\mathbb{R}^{2} ,N)$ with a compact orbit and
$C^{(1,1)}-$close to $\varphi.$ Moreover, the perturbation is
supported on a neighbourhood of a closed orbit of $Y$ where $Y$ is
one of the infinitesimal generating vector fields of $\varphi.$
\end{theorem}

To use the closing lemma of Pugh, we should take care for some
technical problems which arises when we are dealing with
 actions. In fact the idea is as
following: Whenever, we have a dense cylinder we choose a closed
orbit of one of the infintesimal generating vector fields and take
an adequate system of coordinates around this closed orbit. An small
box in our new coordinates will serve as a flow box of the closing
lemma of Pugh. All the orbits passing through this box are two
dimensional. Like in the closing lemma for flows, we have a
transversal section, which is a ring in our case. We construct this
ring foliated by closed orbits of $Y$. Now, we should take care
about the returns of the dense orbit in our neighbourhood. More
precisely, in the following lemma we show that the dense orbit
returns and intersects the transversal section in closed orbits. So,
we really can carry the proof of the closing lemma  for flows to our
case.

The following lemma was announced in \cite{RRW} for non-singular
actions. Here we show that it is true for general actions. In the
following lemma $X$ and $Y$ are two generating infinitesimal vector
fields for $\phi \in A^{1} ( \mathbb{R}^{2} ,M).$

\begin{lemma} \label{le.retorns}
 Let $\mathcal{O}_p$ be a dense cylindrical orbit of $\phi \in A^{1} (
\mathbb{R}^{2} ,M)$ and $c$ (homeomorphic to $S^1$) be the periodic
orbit of $Y$ passing through $p.$ Then, for any neighbourhood $U(c)$
of $c$ there exist an unbounded sequence $t_i \in \RR$ such that
$X^{t_i}(c) \subset U(c).$
\end{lemma}

\begin{proof}
Let $\mathcal{U}_{\epsilon}$ be  an $\epsilon-$neighbourhood of $c$
in $M$ such that $Y^t(z) \in U(c)$ for any $z \in \mathcal{U} , t
\in [0, 1].$ By density of $\mathcal{O}_p$ there exist $z \in c, t
\in \RR$ such that $X^t(z) \in \mathcal{U}_{\epsilon}.$ It comes out
that $Y^I(X^t(z)) \in U(c)$ where $Y^I(.)$ stands for $\{Y^s (.), s
\in I=[0, 1]\}.$ But by commutativity
$$Y^I(X^t(z)) = X^t(Y^I(z)) = X^t(c) \in U(c).$$ As
$\epsilon$ can be any small number, we conclude that there is a
sequence $t_i \ra \infty$ such that $X^{t_i} (c) \in U(c).$
\end{proof}

%

 \noindent
 \vspace{2pt}
 \textbf{Infinitesimal generators adapted to a $S^1 \times \RR$-orbit.}
 Let $\mathcal{O}_{p}$ be a cylindrical orbit of
 $\varphi \in A^{1}(\mathbb{R}^{2},N)$ and $\{w_{1},w_{2}\}$ be a set of generators of its
 isotropy group $G_{p}\,.$ Write
   $X =X_{w_{1}}, Y = X_{w_{2}}$ and $Y$ has periodic orbit of period one through $p.$
    Note that if $q\in \mathcal{O}_p,$ then
    the orbit of $Y$ passing through $q$ is periodic of period one too.
   Put a Riemannian metric on $N$ and let $\xi$ be the norm one vector
   field defined in a neighborhood of the  $Y^{I}_{p}$ that is orthogonal
   to the orbits of $\varphi$.

 Let $c$ be the circle orbit of $Y$
   through $p$. We know that all the orbits close to $c$ are cylindrical or torus. For small $\ep > 0$ we define
   the ring $A_{\ep} = \{ Y^{I} ( \xi^t(c)), |t| \leq \ep \}.$  As the action is $\varphi$ is orientable and $\ep$ is small,
    $A_{\ep}$ is diffeomorphic to $S^1 \times (-\ep, \ep).$ We
    parametrize $c$ with $\theta \in [0, 1]$ such that $\frac{\partial}{\partial \theta} = Y| c.$
   We  put a coordinate system $(x, \theta, z)$ in a small neighbourhood of $c$
   diffeomorphic to $S^1 \times (-\ep, \ep) \times (-1,1)$
   such that:
   $$X = \frac{\partial}{ \partial x},  Y = \frac{\partial}{\partial \theta}, \xi = \frac{\partial}{\partial z}.$$
In this new coordinates system the (pieces of) orbits of $\varphi$
inside such neighbourhood are $z
   =$ constant.

 \subsection{Proof of the main theorem}

 Let $\mathcal{U}$ be a $C^{(1,1)}$ neighborhood of $\varphi$ such
that every action in $\mathcal{U}$ is transitive.

We will prove that $\varphi$ is defined by a flow. It is obvious (by
continuity) that if the dense orbit of a transitive action is
one-dimensional then the action is defined by a flow (The two
infinesimal vestor fieds are linearly dependent.) In what follows we
will show that a robustly transitive action can not have a dense
cylinder. So, we conclude that in fact $\varphi$ is a robustly
transitive flow and by a result of Doering \cite{Do87}, it comes out
that $\varphi$ is an Anosov flow.

Suppose that $\varphi$ has a dense cylinder $\mathcal{O}_p$. Let
$c$ be the periodic orbit (homeomorphic to $S^1$) through $p$ and
$A_{\ep}$ the ring defined previously. Recall that $\{z=0\} \cap
A_{\ep} = c$ and all $\{z= t\} \cap A_{\ep} , |t| \leq \ep$ are
periodic orbits of the generating vector fireld $Y.$ By lemma
\ref{seila2} all two dimensional orbits are either cylindrical or
homeomophic to torus.  Before continuing the argument we state an
algebraic topological lemma. \footnote{The authors would like to
thank C. Biasi for usefull comments and a proof on this lemma.
Later, we find out that a similar lemma was proved in \cite{RR70}.
So, we omit the similar proof.}
\begin{lemma} \label{l.topol}
 Let $N$ be a three dimensional compact orientable manifold. There
 exists $k \in \NN$ such that if $T_1, T_2, \cdots, T_k$ are
 submanifolds homeomorphic to torus $\mathbb{T}^2$, then they form
 the boundary of a three dimensional submanifold of $N.$
\end{lemma}
  Having in mind the above lemma, we conclude that if $A$
has $k$ compact orbits then there can not exist any dense orbit.
Indeed, any dense two dimensional submanifold should intersect one
of these $k$ tori.

So, let  $\mathcal{O}_p$ be a dense orbit and for small $\ep$ all
the orbits passing through $A_{\ep}$ be cylindrical. In fact, if
there does not exist such an $\ep$ we conclude that there are more
than $k$ torus and using the above lemma we contradict the
 denseness of $\mathcal{O}_p.$ For $0
  \leq i \leq k-1$ Let

 $$
  A_i := \big\{ \frac{\ep i}{k} < z < \frac{(i+1) \ep}{k} \big \}.
$$
By the denseness of $\mathcal{O}_p$ and lemma \ref{le.retorns}
there exists a return time $\bar{t}$ such that $X^{\bar{t}}(c) \in
\{ |z| < \frac{\ep}{k}\}.$ As $X^{\bar{t}} (p) \in \{ |z| <
\frac{\ep}{k}\}$ we project $X^{\bar{t}} (p)$ along the orbit of
$X$ and find out $t$ such that $X^{t}(p) \in A_0.$ By definition,
$A_{\ep}$ is foliated by the orbits of $Y$ and by the
commutativity of $X$ and $Y$ one concludes that $X^t (c) \in A_0$
for some $t.$ Indeed,
$$X^t (Y^s(p)) = Y^{s} (X^t (p)) \in A_0 \quad \text{for all} \quad
s \in [0, 1]$$ which means that $X^t(c) \in A_0.$

Now, we use the closing lemma and perturb $\varphi$ inside  $\{
|z| < \frac{\ep}{k} \}$ and find a new action $\varphi_1$ and
$C^{(1,1)}$-close to $\varphi$ with a compact orbit. As
$\varphi_1$ is also transitive, it has a dense orbit which we
claim it is of cylindrical type. To see this remember that our
perturbation is supported on $\{|z| < \frac{\ep}{k}   \}$
  and consequently the orbits passing through $A_i ,
i > 0$ remains cylindrical. So, the dense orbit of $\varphi_1$
which necessarily intersect $\{ \frac{1}{k} <  z < \frac{2}{k}\}$
is cylindrical.  Pertubing again by the closing lemma we obtain
another invariant torus and by induction we find $\varphi_k \in
A^{2}(\mathbb{R}^{2},N)$ with $k$ compact leaves which  by lemma
\ref{l.topol} form the frontier of a compact three manifold with
boundary inside $N$ and so, no dense orbit can exist which gives a
contradiction.

\end{document}